\documentclass[12pt,reqno]{amsart}

\usepackage[centertags]{amsmath}
\usepackage{hyperref}
\usepackage{amsfonts}
\usepackage{amssymb}
\usepackage{amsthm}
\usepackage{newlfont}
\usepackage{amscd}
\usepackage{amsmath,amscd}
\usepackage{enumerate}
\usepackage{color}

\newtheorem{theo}{Theorem}[section]
\newtheorem{prop}[theo]{Proposition}
\newtheorem{lem}[theo]{Lemma}
\newtheorem{cor}[theo]{Corollary}
\newtheorem{defi}[theo]{Definition}
\theoremstyle{remark}
\newtheorem{rem}[theo]{Remark}
\newtheorem{ex}[theo]{Example}

\def \Br {{\mathrm{Br}}}

\def \Ga {{\Gamma}}
\def \R {{\mathbb{R}}}
\def \Pic {{\rm {Pic}}}

\def \Gal {{\mathrm{Gal}}}
\def \Ker {{\mathrm{Ker}}}

\def \AA{{\mathbf{A}}}
\def \A{{\mathbb A}}
\def \P{{\mathbb P}}
\def \Spec {{\mathrm{Spec}}}
\def \dim {{\mathrm{dim}}}
\def \Hom {{\mathrm {Hom}}}

\def \Pic {{\mathrm {Pic}}}
\def \GL {{\mathrm{GL}}}
\def \SL {{\mathrm{SL}}}

\def\ov{\overline}

\def \Z {{\mathbb Z}}
\def \Q {{\mathbb Q}}

\def \Id {{\mathrm{Id}}}

\def\G{{\mathbb G}}

\def\lra{\longrightarrow}

\def\res{{\mathrm{res}}}

\def\Ga{\Gamma}

\def\e{\varepsilon}
\def\et{\mathrm{\acute et}}

\newcommand{\bthe}{\begin{theo}}
\newcommand{\ble}{\begin{lem}}
\newcommand{\bpr}{\begin{prop}}
\newcommand{\bco}{\begin{cor}}
\newcommand{\bde}{\begin{defi}}
\newcommand{\ethe}{\end{theo}}
\newcommand{\ele}{\end{lem}}
\newcommand{\epr}{\end{prop}}
\newcommand{\eco}{\end{cor}}
\newcommand{\ede}{\end{defi}}

\title{Pathologies of the Brauer--Manin obstruction}

\date{\today}
\begin{document}
\author{Jean-Louis Colliot-Th\'el\`ene}
\address{CNRS, UMR 8628, Math\'ematiques, B\^atiment 425, Universit\'e Paris-Sud,
F-91405 Orsay, France}
\email{jlct@math.u-psud.fr}

\author{Ambrus P\'al}
\address{Department of Mathematics, South Kensington Campus,
Imperial College London, SW7 2BZ England, U.K.}
\email{a.pal@imperial.ac.uk}

\author{Alexei N. Skorobogatov}
\address{Department of Mathematics, South Kensington Campus,
Imperial College London, SW7 2BZ England, U.K. - and - 
Institute for the Information Transmission Problems,
Russian Academy of Sciences, 19 Bolshoi Karetnyi, Moscow, 127994
Russia}
\email{a.skorobogatov@imperial.ac.uk}

\baselineskip=15pt
\maketitle

\section{Introduction}

It had long been suspected that varieties $X$ over a number field $k$
without $k$-points but with a non-empty Brauer--Manin set
$X(\AA_k)^{\Br}$ are fairly common.
The first examples were found in \cite{Sk1, BSk}.
An earlier, conditional example is given  in \cite{SW}. One should also expect
that there are many varieties $X$ without $k$-points for which the
\'etale Brauer--Manin set $X(\AA_k)^{\et,\Br}\subset X(\AA_k)^{\Br}$ 
is non-empty. (We refer to \cite{P1} or \cite{Sk2} 
for the definition of these
subsets of the space $X(\AA_{k})$ of ad\`eles of $X$.)
Different methods to construct such varieties have been found
recently. In \cite{P1} Poonen 
constructs a threefold $X$ with a surjective morphism 
to a curve $C$ that has exactly one $k$-point $P$ and
the fibre $X_P$ has points everywhere locally 
but not globally. In Poonen's example $X_P$
is a smooth Ch\^atelet surface. The trick with a curve with just
one rational point was also used in
\cite{HaSk} where the fibres of $X\to C$ are curves of high genus and 
$X_P$ is a singular curve which geometrically is a union of 
projective lines.
In retrospect one could note that the examples in \cite{Sk1} and \cite{BSk}
are families of genus 1 curves parameterised by elliptic curves
of Mordell--Weil rank~0.

In this paper we propose more flexible methods to construct such examples.
We show that the varieties $X$ such that $X(k)=\emptyset$ and
$X(\AA_k)^{\et,\Br}\not=\emptyset$ include the following:

a conic bundle surface $X\to E$ over a real quadratic field $k$,
where $E$ is an elliptic curve such that $E(k)=\{0\}$,
see Section \ref{dim2};

a threefold over an arbitrary real number field $k\subset\R$, 
which is a family $X \to C$ of 2-dimensional quadrics parameterised by a
curve $C$ with exactly one $k$-point (one can choose
$C$ to be an elliptic curve), see Section \ref{dim3};

a threefold over an arbitrary number field $k$,
which is a family $X \to C$ 
of geometrically rational surfaces parameterised by a
curve $C$ with exactly one $k$-point, the fibre above which
is singular, see Section \ref{def}.

In the first and second examples, 
in contrast to those previously known, 
the smooth fibres satisfy the Hasse principle and weak approximation.
To put this into a historical perspective let us note that soon after
Manin \cite{Manin} introduced the obstruction now bearing his name,
Iskovskikh \cite{Isk} 
constructed a counterexample to the Hasse principle 
on a conic bundle over the projective line over $\Q$. 
His intention was, as he
pointed out to one of us, to give a counterexample to the Hasse principle
that could not be explained by the Brauer--Manin obstruction.
It is well known nowadays that Iskovskikh's counterexample 
can be explained by the Brauer--Manin obstruction,
and conjecturally the same should be 
true for all counterexamples to the Hasse principle on 
geometrically rational surfaces, see \cite{angers, CTCoSa}.

The examples
we construct in this paper show that this is no longer the case
for conic or quadric bundles over curves of genus at least 1.

In a nutshell, the idea is this. Let $k$ be a number field.
Following Poonen we use a base variety $B$ such that $B(k)=\{P\}$. 
By a continuous deformation of the ad\`ele attached to $P$ at
an archimedean component
we see that $B(k)$ is not dense in $B(\AA_k)^{\Br}$.
Density may also fail due to places of $k$ that
need not be archimedean.
Suppose $B$ contains an irreducible singular conic
$S$ so that $P=S_{\mathrm{sing}}$. If a place $v$ of $k$ 
splits in the quadratic extension given by the discriminant of 
the binary quadratic form that defines $S$,
then $B\times_k k_v$ contains two copies of $\P^1_{k_v}$ meeting at $P$.
Since $\Br(\P^1_{k_v})=\Br(k_v)$,
we can modify the ad\`ele of $P$ at $v$ while staying inside $B(\AA_k)^{\Br}$.
However, the $k$-point $P$ cannot be moved in $B$, 
so $B(k)$ is not dense in $B(\AA_k)^{\Br}$.

Next, one constructs a surjective morphism $X\to B$ for which the fibre $X_P$
has local points in all but one or two completions of $k$, and ensures that
$X$ has $k_v$-points for missing places $v$ such the resulting adelic point
of $X$ projects to $B(\AA_k)^{\Br}$. Now, if the natural map $\Br(B)\to\Br(X)$
is surjective we have found an adelic point in $X(\AA_k)^{\Br}$. But since
$X(k)\subset X_P$ we have $X(k)=\emptyset$. With more work one can find
examples such that $X(\AA_k)^{\et,\Br}$ is non-empty, too.

In this paper we have nothing to say about the important open question
whether the implication
$$X(\AA_k)^{\et,\Br}\not=\emptyset\quad\Rightarrow \quad X(k)\not=\emptyset$$
holds if $X$ is a surface with finite geometric fundamental group,
e.g. a K3 surface or an Enriques surface. 

The paper is organised as follows. After some preparations in Section
\ref{prep} we realise the aforementioned programme
for threefolds in Section \ref{3folds}. 
Making it work for surfaces requires rather more effort.
For this purpose in Section \ref{ec} we establish some 
Bashmakov-style properties of elliptic curves 
with a large Galois image on torsion points.
These properties are used in the proof of our main result 
in the case of surfaces in Section \ref{surfaces}. Some general observations
on the Brauer--Manin set are collected in Section \ref{varia}.

The authors are grateful to Nathan Jones and Chris Wuthrich 
for their helpful advice on Serre curves. We used {\tt sage}
in our calculations with elliptic curves. This work started 
in April 2013 when the authors
were guests of the Hausdorff Institut f\"ur Mathematik (Bonn) during
the special programme ``Arithmetic and geometry".

\section{Brauer groups and torsors on quadric bundles} \label{prep}

For the convenience of the reader we recall the following well known lemma.

\ble \label{n1}
Let $k$ be a field of characteristic zero. Let $X$ be a smooth
projective quadric over $k$ of dimension at least $1$.
Then the natural map $\Br(k)\to\Br(X)$ is surjective.
If $\dim(X)\geq 3$, then this map is an isomorphism.
\ele
{\em Proof.} Let $\bar k$ be an algebraic closure of $k$,
and let $\Ga_k=\Gal(\bar k/k)$. 
For any smooth, projective and geometrically integral variety $X$
over $k$ there is a well known exact sequence
$$ 0 \to\Pic(X) \to\Pic({\ov X})^{\Ga_k} \to \Br(k) \to 
\Ker[\Br(X) \to \Br(\ov X)] \to H^1(k,\Pic(\ov X)),$$
where $\ov X=X\times_k\bar k$. If $X$ is a quadric of
dimension at least $1$, then
$\Pic(\ov X)$ is a permutation $\Ga_k$-module
and $\Br(\ov X)=0$ (since $\ov X$ is
rational and the Brauer group is a birational invariant
of smooth projective varieties).
By Shapiro's lemma we have $ H^1(k,\Pic(\ov X))=0$,
so the exact sequence implies the surjectivity of the map $\Br(k)\to\Br(X)$.
When $\dim(X)\geq 3$, the map $\Pic(X)\to\Pic(\ov X)$ is an isomorphism,
because both groups are generated by the hyperplane section class,
so in this case $\Br(k)\to\Br(X)$ is an isomorphism. $\Box$

\medskip

In this paper a {\em quadric bundle}
is a surjective flat morphism $f : X \to B$ of smooth, projective,
geometrically integral varieties over a field $k$,
the generic fibre of 
which is a smooth quadric of dimension at least $1$,
and all geometric fibres are reduced.
  
We denote by $k(B)$ the function field of $B$, and by $X_{k(B)}$
the generic fibre of $f : X \to B$. If $\dim(X_{k(B)})=1$,
then $f : X \to B$ is called a {\em conic bundle}.

The following proposition is essentially well known, at least
when $B=\P^1_k$, see \cite[Cor. 3.2]{Sk90},
\cite[Thm. 2.2.1, Thm. 2.3.1]{CTSD94}, \cite[Prop. 2.1]{CTSMF}.

\bpr\label{Br-const}
Let $f : X \to B$ be a quadric bundle 
over a field $k$ of characteristic zero. 
In each of the following cases the map $f^{*}:\Br(B)\to\Br(X)$ is surjective.

$(i)$ $\dim(X_{k(B)})=1$ 
and there is a point $P\in B$ of codimension $1$ such that 
for each point $Q\not=P$ of codimension $1$ in $B$ the fibre $X_Q$
contains a geometrically integral component of multiplicity $1$;

$(ii)$ $\dim(X_{k(B)})=2$ and for each
point $Q\in B$ of codimension $1$ the fibre $X_Q$
contains a geometrically integral component of multiplicity $1$;

$(iii)$ $\dim(X_{k(B)})\geq 3$.
\epr
{\em Proof.} $(i)$ Let $\gamma \in \Br(k(B))$ be the class of the conic
$X_{k(B)}$. Since $\gamma$ is in the kernel of the natural map
$f^*:\Br (k(B)) \to \Br(X_{k(B)})$, the assumption of $(i)$ 
implies that 
the residue $\res_Q(\gamma)\in H^1(k(Q),\Q/\Z)$ is zero if $Q\not=P$.
Take any $\alpha \in \Br(X)$.
By Lemma \ref{n1} the map $f^*:\Br (k(B)) \to \Br(X_{k(B)})$ is surjective,
so the image of $\alpha$ in $\Br(X_{k(B)})$ comes from some 
$\beta \in \Br(k(B))$. Again, by the assumption of $(i)$
we have $\res_Q(\beta)=0$ if $Q\not=P$. Moreover,
$\res_P(\beta)=0$ or $\res_P(\beta)=\res_P(\gamma)$. By the purity theorem
for the Brauer group \cite[III, Thm. 6.1, p. 134]{GB}
we conclude that $\beta\in\Br(B)$ or $\beta-\gamma\in\Br(B)$. 
Since $f^*(\beta)=f^*(\beta-\gamma)=\alpha$ in $\Br(X_{k(B)})$,
and the natural map $\Br(X)\to \Br(X_{k(B)})$ is injective, 
we have proved $(i)$. 

The proof of $(ii)$ uses Lemma \ref{n1}
and the arguments from the proof of $(i)$.

In case $(iii)$ it is well known that each fibre of $f$ at a point 
of $B$ of codimension $1$ contains
a geometrically integral component of multiplicity $1$. Then 
$(iii)$ follows from the last statement of Lemma \ref{n1}. $\Box$

\bpr  \label{et-cov}
Let $f : X \to B$ be a quadric bundle over 
a field $k$ of characteristic zero.
 Then any torsor
$X' \to X$ of a finite $k$-group scheme $G$ is the inverse image under $f$
of a torsor $B' \to B$ of $G$.
\epr
{\em Proof}. 
By our definition of quadric bundles, the morphism $f$ 
is flat and  all its geometric  fibres are connected and reduced. 
The generic geometric fibre of $f$ is   simply connected.
By \cite[X, Cor. 2.4]{SGA1} this implies  that  each geometric fibre 
of such a fibration is simply connected.
The result then follows from \cite[IX, Cor. 6.8]{SGA1}. $\Box$

\section{Threefolds} \label{3folds}

\subsection{Example based on real deformation} \label{dim3}

Let $k$ be a number field with a real place.
We fix a real place $v$, so we can think of $k$ as a subfield of $k_v=\R$.

Let $C$ be a smooth, projective, geometrically integral curve over 
$k$ such that $C(k)$ consists of just one point, $C(k)=\{P\}$.
By \cite{P2} such a curve exists for any number field $k$,
and by \cite[Thm. 1.1]{MR} we can take $C$ to be an elliptic curve
over $k$.
Let $\Pi \subset C(\R)$ be an open interval containing $P$. 
Let $f:C \to \P^1_k$ be a surjective morphism that is
unramified at $P$.
Choose a coordinate function $t$ on $\A^1_k=\P^1_k\setminus f(P)$
such that $f$ is unramified above $t=0$. We have $f(P)=\infty$. 
Take any $a>0$ in $k$ such that
$a$ is an interior point of the interval $f(\Pi)$ 
and $f$ is unramified above $t=a$.

Let $w$ be a finite place of $k$. There exists a quadratic form
$Q(x_0,x_1,x_2)$ of rank $3$ that represents zero in all
completions of $k$ other than $k_v$ and $k_w$, but not in 
$k_v$ or $k_w$. We can assume that $Q$ is positive definite
over $k_v=\R$. Choose $n \in k$ with $n>0$ in $k_{v}$
and $-nQ(1,0,0)\in k_w^{*2}$.
Let $Y_1\subset\P^3_k\times\A^1_k$ be
given by $Q(x_0,x_1,x_2)+nt(t-a)x_3^2=0$, and let
$Y_2\subset\P^3_k\times\A^1_k$ be
given by $Q(X_0,X_1,X_2)+n(1-aT)X_3^2=0$. We glue
$Y_1$ and $Y_2$ by identifying
$T=t^{-1}$, $X_3=tx_3$, and $X_i=x_i$ for $i=0,1,2$.
This produces a quadric bundle $Y\to \P^1_k$ with exactly
two degenerate fibres (over $t=a$ and $t=0$), 
each given by the quadratic form $Q(x_0,x_1,x_2)$ of rank 3. 
Define $X=Y\times_{\P^1_k} C$. 
This is a quadric bundle $X\to C$ with
geometrically integral fibres.

For example, if $k=\Q$, we can take $k_w=\Q_2$ and consider $Y$ defined by
$$x_0^2+x_1^2+x_2^2+7t(t-a)x_3^2=0.$$

\bpr \label{real}
In the above notation we have $X(\AA_k)^{\et,\Br}\not=\emptyset$ 
and $X(k)=\emptyset$.
\epr
{\em Proof}. 
Since $C(k)=\{P\}$ we have $X(k)\subset X_P$. The fibre $X_P$
is the smooth quadric $Q(x_0,x_1,x_2)+nx_3^2=0$.
This quadratic form is positive definite
thus $X_P$ has no points in $k_v=\R$ and so $X(k)=\emptyset$.
By assumption $X_P$ has local points in all completions 
of $k$ other than $k_v$ and $k_w$.
The condition $-nQ(1,0,0)\in k_w^{*2}$ implies that $X_P$
contains $k_w$-points, so $X_P$ has local points in all completions
of $k$ but one.
Choose $N_u\in X_P(k_u)$ for each place $u\not=v$.
Consider a small real $\e>0$ such that $a-\e \in f(\Pi)$ and $\e<a$.
Let $M\in\Pi$ be such that $f(M)=a-\e$. Then
the smooth real fibre $X_M$ is given by an indefinite quadratic form
and so $X_M(k_v)\not=\emptyset$. Choose any $N_v \in X_M(k_v)$.
We now have an adelic point $(N_u)$, where we allow $u=v$.

We claim that $(N_u)\in X(\AA_k)^{\et,\Br}$.

Let $G$ be a finite $k$-group scheme. 
Proposition \ref{et-cov} implies that any torsor $X'/X$ of
$G$ comes from a torsor $C'/C$ of $G$, 
in the sense that $X\times_C C'\to X$ 
and $X'\to X$ are isomorphic as $X$-torsors with the structure group $G$.
Let  $\sigma \in Z^1(k,G)$ be a 1-cocycle defining the $k$-torsor 
which is the fibre of $C' \to C$ at $P$.
Twisting $X'/X$ and $C'/C$ by $\sigma$ and replacing the group $G$
by the twisted group $G^{\sigma}$ and changing notation,
we can assume that $C'$ contains a $k$-point $P'$
that maps to $P$ in $C$. The irreducible 
component $C''$ of   $C'$ that contains $P'$
is a geometrically integral curve over $k$. 
Let $X'' \subset X'$ denote the inverse
image of $C''$ in $X'$. The fibres of the morphism $X\to C$ are 
geometrically integral, hence such are also the fibres of $X' \to C'$ 
and $X'' \to C''$. Thus $X''$ is a geometrically integral variety
over $k$.
 
There are natural isomorphisms $X''_{P' } \cong X'_{P'} \cong X_P$, 
so we can define $N'_u\in X''(k_u)$ as the point 
that maps to $N_u\in  X(k_u)$ for each $u\not=v$. The map 
$C'' \to C$ is finite and \'etale. 
The image of $C''(\R)$ in $C(\R)$ is thus closed and open.
The image of the connected component of $P' \in C''(\R)$ 
is the whole connected component of $P \in C(\R)$, hence contains $\Pi$.
The inverse image of the interval $\Pi$ in $C''(\R)$
is a disjoint union of intervals, one of which contains $P'$
and maps bijectively onto~$\Pi$.
Let us call this interval $\Pi'$. 
Let $M'$ be the unique point of $\Pi'$ over $M$. 
Let $N'_v\in X''_{M'}(\R)$ be the point that maps to 
$N_v\in X_M(\R)$. 
Thus the adelic point $(N'_u)\in X''(\AA_k)\subset X'(\AA_k)$ 
projects to the adelic point $(N_u)\in X(\AA_k)$.

By the definition of the \'etale Brauer--Manin obstruction,
to prove that $(N_u)$ is contained in $X(\AA_k)^{\et,\Br}$ 
it suffices to show that $(N'_u)$ is orthogonal to $\Br(X')$.
For this it is enough to show that $(N'_u)$ is orthogonal to $\Br(X'')$.
By Proposition \ref{Br-const} $(ii)$ applied to  $X'' \to C''$
we know that the natural map
$\Br(C'') \to \Br(X'')$ is surjective. Thus it is enough to show that
the ad\`ele on $C''$ such that its $u$-adic component is $P'$ when $u\not=v$
and and its $v$-component is $M'$,
is orthogonal to $\Br(C'')$. The real point $M'$ is path-connected to $P'$,
so this ad\`ele is in the connected component of  
the diagonal image of the $k$-point $P'$ in $C''(\AA_k)$. But
the latter ad\`ele is certainly in $C''(\AA_k)^\Br$,
and the proposition follows. $\Box$

\begin{rem} \label{r2}
(1) Our method gives
simple examples of threefolds with points everywhere locally
but not globally and no Brauer--Manin obstruction. An even simpler
proof is available in the case of fibrations 
into quadrics of dimension at least $3$ over a curve.

(2) By a theorem of Wittenberg \cite[Thm. 1.3]{W} the variety $X$ has 
a $0$-cycle of degree $1$ over $k$, that is, there exist
field extensions $k_1,\ldots,k_r$ of $k$ whose degrees
have no common factor such that $X(k_i)\not=\emptyset$
for $i=1,\ldots, r$. Although \cite[Thm. 1.3]{W} requires
the finiteness of the Shafarevich--Tate group of the Jacobian of $C$, 
Wittenberg pointed out that in the proof of his theorem 
this assumption is only used to ensure the existence
of a suitable $0$-cycle of degree $1$ on $C$. In our case 
such a $0$-cycle is directly provided by the $k$-point $P$, so
the assumption on the  Shafarevich--Tate group is not needed. For more
details see Remark \ref{r3} below.
\end{rem}

\subsection{Examples based on deformation along a rational curve
defined over a completion of $k$} \label{def}

\ble \label{lle}
Let $k$ be a number field. There exists a smooth, projective,
geometrically integral surface $B$ over $k$ with the following properties:

$B$ contains a curve $S$ isomorphic to an irreducible singular 
projective conic;

the singular point of $S$ is the unique $k$-point of $B$;

there is a surjective morphism $\pi:B\to\P^1_k$ with smooth and geometrically
integral generic fibre such that $\pi(S)=\P^1_k$.
\ele
{\em Proof.} According to \cite{P2} there is
a smooth, projective and geometrically integral curve $C$ over $k$
with exactly one $k$-point, $C(k)=\{O\}$. Moreover,
by \cite[Thm 1.1]{MR} there is an elliptic curve over $k$ with this property.
Let $f:Z\to C$ be any conic bundle such that
the fibre $Z_O$ over $O\in C(k)$ is an irreducible singular conic. 
The singular point of $S=Z_O$ is then the unique $k$-point of $Z$.

There is a closed embedding $Z\subset\P^n_k$ for some 
$n\geq 1$. By the Bertini theorem \cite[II.8.18, III.7.9]{Hartshorne}, 
there exists a hyperplane $H_1\subset \P^n_k$
such that $Z\cap H_1$ is a smooth and geometrically integral curve.
This implies that $S$ is not a subset of $H_1$.

Let $d$ be the degree of $Z$ in $\P^n_k$. We can find 
a hyperplane $H_2\subset\P^n_k$ such that $Z\cap H_1\cap H_2$ is 
a set of $d$ distinct $\bar k$-points not
in $S$. It follows that no geometric irreducible
component of $S$ is contained in a hyperplane passing through $H_1\cap H_2$.

Let $\tilde\P^n_k$ be the blowing-up of $\P^n_k$ at 
$H_1\cap H_2\simeq\P^{n-2}_k$. The projection from $H_1\cap H_2$
defines a morphism $\tilde\P^n_k \to \P^1_k$. Let
$B$ be the Zariski closure of $Z\setminus (Z\cap H_1\cap H_2)$
in $\tilde\P^n_k$, so that $B$ is the blowing-up of $Z$
in $d$ distinct points. Thus $B$ is a smooth, projective,
geometrically integral surface with a unique $k$-point, 
equipped with a surjective morphism $\pi:B\to\P^1_k$
with smooth and geometrically integral generic fibre.
Moreover, $S$ is contained in $B$, and $\pi(S)=\P^1_k$. $\Box$

\medskip

Let $P$ be the unique $k$-point of  $B$, and let $Q=\pi(P)\in\P_k^1(k)$.
Let $K$ be the quadratic extension of $k$ over which the components
of $S$ are defined. If $w$ is a place of $k$ that splits in $K$,
then the $k_w$-variety $S\times_k k_w$ is the union of 
two projective lines meeting at $P$. 
Let $L_w\subset B\times_k k_w$ be one of these rational curves. 
Since $\pi(L_w)=\P^1_{k_w}$, 
there is a point $N_w\in L_w(k_w)$ such that $\pi(N_w)\not=Q$.

\bpr \label{p-adic}
Let $w_1$ and $w_2$ be places of $k$ that split in $K$, $w_1\not=w_2$.
Let $Y\to\P^1_k$ be a conic bundle satisfying the following conditions:

there exists a closed point $R\in \P^1_k$, $R \neq Q$, 
such that 
the restriction $Y\setminus Y_R \to \P^1_k\setminus R$ is a smooth morphism,
and the fibre of $\pi:B\to\P^1_k$ at $R$ is smooth;

the fibre $Y_Q$ is a smooth conic that has $k_v$-points 
for all completions of $k$
except $w_1$ and $w_2$, in particular $Y_Q(k)=\emptyset$;

$Y_{\pi(N_w)}(k_{w})\not=\emptyset$ for $w=w_1, w_2$. 

\noindent Then for the smooth  threefold $X=Y\times_{\P^1_k}B$ we have
$X(\AA_k)^{\et,\Br}\not=\emptyset$ and $X(k)=\emptyset$.
\epr
{\em Proof}. Let $p:X\to B$ be the natural projection.
Since $B(k)=\{P\}$, $\pi(P)=Q$ and $Y_Q(k)=\emptyset$,
we see that $X(k)=\emptyset$.

The fibre $X_P$ is naturally isomorphic to $Y_Q$.
For a place $v$ such that $v\not=w_1$, $v\not=w_2$ choose $M_v\in X_P(k_v)$.
For $w=w_1,w_2$, choose a $k_{w}$-point $M_{w}$ in the $k_{w}$-fibre $X_{N_{w}}$
(which is isomorphic to $Y_{\pi(N_{w})}$). We claim that
$(M_v)\in X(\AA_k)^{\et,\Br}$.

One easily checks that the projection map $X \to B$ 
is a quadric bundle
as defined in this paper: both $X$ and $B$ are smooth and projective 
over $k$, the generic fibre is a geometrically
integral conic, the morphism $X \to B$ is flat and all its geometric 
fibres are reduced.
Let $G$ be a finite $k$-group scheme.
By Proposition \ref{et-cov} every torsor $X'\to X$ of $G$
is the pullback of a torsor $B'\to B$ of $G$.
After a twist by a $k$-torsor of $G$, as detailed in the proof of Proposition
\ref{real}, we may assume that $B'$ has a $k$-point $P'$ over $P$. 
Let $B''$ be the irreducible component of $B'$ that contains $P'$.
Then $B''$ is geometrically integral. Let
$X''$ be the inverse image of $B''$ under the map $X'\to B'$.
The $k$-variety $X''$ is geometrically integral, and
$X''\to B''$ is a conic bundle.

For all $v\not=w_1,w_2$ the fibre $X''_{P'}$
contains a $k_v$-point $M'_v$ that maps to $M_v$ in $X_P$.
Since $B''\to B$ is finite and \'etale, for $w=w_1,w_2$
the inverse image of $L_{w}$
in $B''$ contains a rational $k_{w}$-curve $L'_{w}$ through $P'$
that maps isomorphically to $L_{w}$. Let $N'_{w}$ be the $k_{w}$-point
of $L'_{w}$ that maps to $N_{w}$. Then there is a $k_{w}$-point
$M'_{w}$ in $X''_{N'_{w}}$ that maps to $M_{w}\in X_{N_{w}}(k_{w})$.
To prove our claim it is enough to show that
the adelic point $(M'_v)$ in $X''$ is orthogonal to $\Br(X'')$.
 
Let $\bar x\in R$ be a $\bar k$-point. Since the singular loci of the 
morphisms $\ov Y\to\P^1_{\bar k}$ and $\pi:\ov B\to\P^1_{\bar k}$ are disjoint,
the fibre $\ov B_{\bar x}=\pi^{-1}(\bar x)$ 
is a smooth integral curve over $\bar k$.
The fundamental group of $\P^1_{\bar k}$ is trivial. Therefore,
by \cite[Chap. X, Cor. 1.4]{SGA1}, the map
$\pi_1(\ov B_{\bar x})\to \pi_1(\ov B)$ is surjective.
Since $\ov B''$ is integral, this implies that
$B''\times_B B_{\bar x}$ is also integral. As a consequence, 
the inverse image of $R$ in $B''$ is a geometrically integral
$k(R)$-curve $D$. Hence $D$ is integral.

Since the restriction of the conic bundle $X''\to B''$  
to the complement $B''\setminus D$ of 
the integral curve $D$ is a smooth morphism, 
this conic bundle satisfies the condition of Proposition \ref{Br-const} $(i)$.
It follows that the induced map $\Br(B'') \to \Br(X'')$ is surjective.
Thus it remains to prove that the image of $(M'_v)$ in $B''$
is orthogonal to $\Br(B'')$.

This image is the adelic point
such that for all $v\not=w_1,w_2$ the $k_v$-component is $P'$.
For $w=w_1, w_2$ the $k_w$-component is $N'_w$. But $L'_w$ is
a projective line over $k_w$ containing both $P'$ and $N'_w$.
The natural map $\Br(k_w)\to\Br(L'_w)$ is an isomorphism.
Thus when pairing the image of $(M'_v)$ in $B''$
with $\Br(B'')$ we may replace $(M'_v)$ by the diagonal
ad\`ele $(P')$, which by the reciprocity law is orthogonal to $\Br(B'')$.
$\Box$

\begin{rem} Consider a surface $B$
with a morphism $\pi:B\to\P^1_k$ as in Lemma \ref{lle}.
Let $Y \to \P^1_k$ be a quadric bundle all fibres of which
are of dimension $d\geq 2$ and contain a geometrically integral component
of multiplicity one (which is automatic if $d \geq 3$). 
Assume that the singular loci of $B\to\P^1_k$ and $Y\to\P^1_k$
do not intersect. Let $X=Y \times_{\P^1_k}B$. Suppose that the fibre
$Y_Q$, where $Q=\pi(P)\in\P^1(k)$, is a smooth quadric such that 
$Y_Q(k_v)\not=\emptyset$
for all $v\not=w$ but $Y_Q(k_w)=\emptyset$ (such quadrics exist
in dimension 2 and higher, but not in dimension 1).
Assume that $L_w$ has a non-empty intersection with the image
of $X(k_w) \to B(k_w)$. 
In view of Proposition \ref{Br-const} $(ii)$, $(iii)$, an argument
similar but shorter than the one above shows
that $X(k)=\emptyset$ and $X(\AA_k)^{\et,\Br}\not=\emptyset$. 
In this example $\dim(X)=2+d\geq 4$.
\end{rem}

\begin{ex} One can construct  
a threefold over $k=\Q$ as in Proposition \ref{p-adic}
as follows. Let $E$ be the elliptic curve $y^2=x^3-5$.
Then $E(\Q)=\{O\}$, where $O$ is the point at infinity.
Let $(r:s:u:v:t)$ be homogeneous coordinates in $\P^4_\Q$.
The first of the equations
\begin{equation}
xy t^2=u^2+v^2,\quad\quad
r^2-5s^2-17t^2-u^2=0 \label{eqq}
\end{equation}
defines a closed subset of $(E\setminus\{O\})\times\P^2_\Q$ 
which extends to a conic bundle surface $B\to E$.
The fibre $B_O$ over $O$ is the singular conic $S$
with equation $u^2+v^2=0$. The unique $\Q$-point $P$ of $B$
is the singular point of $S$.
The morphism $\pi:B\to\P^1_\Q$ given
by the projection to $(u:t)$ satisfies the conclusions of Lemma \ref{lle}.
The second equation of (\ref{eqq}) defines a smooth quadric 
$Q\subset \P^3_\Q$.
Let $Y\to Q$ be the blowing-up of the closed point of $Q$
given by $u=t=0$. The projection
via the coordinates $(u:t)$ is a morphism
$Y\to\P^1_\Q$ which makes $Y$ a conic bundle as in Proposition \ref{p-adic}. 
Let $X=Y \times_{\P^1_k}B$.
The fibre $X_P$ over $P$ is the conic $r^2-5s^2-17t^2=0$ over $\Q$, so
for the places $w_1$ and $w_2$ one takes 
the primes $5$ and $17$. For $p=5, 17$ choose  
a suitable point $N_p\in S(\Q_p)$ with 
$v=1$ and $u=\alpha_{p}\in\Q_p$ such that $\alpha_{p}^2=-1$.

One can give a different proof of the non-emptyness of the set $X(\AA_\Q)^{\et,\Br}$
using the method of \cite{HaSk}. (By Proposition \ref{5.1} $(ii)$ below
it does not matter which birational model is used for this.)
Let $K=\Q(\sqrt{-1})$. Consider the fibre $X_O$ of the composed 
morphism $X\to B\to E$ over $O\in E(\Q)$. 
The singular surface $X_O$ is fibred into conics over
the singular conic $S$; the inverse image of the singular point 
$P\in S$ is a smooth conic $X_P\subset X_O$.
Thus $X_O\times_\Q K$ is the union of two geometrically
irreducible components permuted by $\Gal(K/\Q)$
that intersect transversally in $X_P$, each of them isomorphic to
$Y_K=Y\times_\Q K$. 
Since $5$ and $17$ split in $K$ and 
the components of $X_O\times_\Q K$ have
$K$-points, we see that $X_O(\AA_\Q)\not=\emptyset$. 
We claim that 
$$X_O(\AA_\Q)\subset X(\AA_\Q)^{\et,\Br}.$$
Indeed, let $G$ be a finite $k$-group scheme. 
The generic fibre of $X\to E$ is a geometrically integral, smooth
and geometrically rational surface, so it is geometrically simply connected.
One checks that the morphism $X\to E$ is flat and all its geometric 
fibres are connected and reduced. Now
\cite[X, Cor. 2.4]{SGA1} implies that each geometric fibre 
of $X\to E$ is simply connected. As in the proof of
Proposition \ref{et-cov} we see from \cite[IX, Cor. 6.8]{SGA1}
that any torsor $X'/X$ of $G$ is the pullback of a torsor $E'/E$ of $G$.
As in the proofs of Propositions \ref{real} and \ref{p-adic} 
it is enough to assume that $E'$ has a 
$\Q$-point $O'$ over $O\in E(\Q)$. Thus a natural isomorphism 
$X'_{O'}\cong X_O$ gives an identification $X'_{O'}(\AA_\Q)\cong X_O(\AA_\Q)$,
so to prove our claim it is enough to show that the natural map
$\Br(\Q)\to\Br(X_O)$ is an isomorphism.

Let $i:X_P\to X_O$ be the closed embedding. Let
$\nu:Y_K\to X_O$ be the normalisation morphism and let $C=\nu^{-1}(X_P)$.
Then $C=X_P\times_\Q K$ is the intersection of the quadric $Q_K=Q\times_\Q K$ 
given by the second equation in (\ref{eqq}) with the plane $u=0$.
The morphism $\nu:C\to X_P$ is the natural projection $X_P\times_\Q K\to X_P$.

The exact sequence of \'etale sheaves on $X_O$
$$0\to\G_{m,X_O}\to \nu_*\G_{m,Y_K}\oplus i_*\G_{m,X_P}\to i_*\nu_*\G_{m,C}
\to 0$$
is similar to the exact sequence (2) in \cite{HaSk}. 
The normalisation morphism $\nu$ and the closed embedding $i$ are
finite morphisms, so $\nu_*$ and $i_*$ are exact functors, hence on
taking cohomology we obtain an exact sequence
$$\Pic(Y_K)\oplus\Pic(X_P)\to\Pic(C)\to \Br(X_O)\to\Br(Y_K)\oplus\Br(X_P)
\to \Br(C).$$
The discriminant of the quadratic form defining $Q$ is not a square
in $K$, hence $\Pic(Q_K)$ is generated by the class of the hyperplane section,
and the natural map $\Br(K)\to \Br(Q_K)$ is an isomorphism.
By the birational invariance of the Brauer group we obtain that
the natural map $\Br(K)\to \Br(Y_K)$ is also an isomorphism.
It is well known that $\Br(X_P)$ is the
quotient of $\Br(\Q)$ by the subgroup generated by 
the class of the conic $X_P$, which is given by the symbol $(5,17)$.
This symbol remains non-zero in $\Br(K)$, hence $C(K)=\emptyset$
and so $\Pic(C)$ is also generated by the class of 
the hyperplane section. Since
the composition of the embedding $C\to Y_K$ with the birational morphism 
$Y_K\to Q_K$ is the natural embedding of $C$ as a plane section of $Q_K$,
we see that the restriction map $\Pic(Y_K)\to\Pic(C)$ is surjective.
Now using the fact that $\Br(C)$ is the quotient of $\Br(K)$ by the subgroup
generated by the symbol $(5,17)$
we easily deduce that $\Br(X_O)=\Br(\Q)$.

To conclude, our example resembles that of Poonen \cite{P1} in that
the fibres of $X\to E$ are birationally equivalent
to intersections of two quadrics in
$\P^4$. However, in our case the fibre above the unique $\Q$-point
is geometrically simply connected and
satisfies $\Br(X_O)=\Br(\Q)$ and $X_O(\AA_\Q)\not=\emptyset$, so
we see that the \'etale Brauer--Manin set of $X$ is non-empty
without studying $X$ any further. Note that 
it is essential that the fibres of $X\to E$ have dimension at least 2.
Indeed, each everywhere locally solvable
geometrically connected and simply connected {\em curve}
over a number field $k$ has a $k$-point \cite[Remark 2.2]{HaSk}.
\end{ex}

\section{Elliptic curves with a large Galois image} \label{ec}

Let $k$ be a field of characteristic zero, with an algebraic closure
$\bar k$. For a field $K$ such that $k\subset K\subset \bar k$ 
we denote by $\Ga_K$ the Galois group $\Gal(\bar k/K)$.
Let $k^{\mathrm{cyc}}\subset \bar k$ be the cyclotomic extension of $k$, 
i.e. the abelian extension of $k$ obtained by adjoining  all roots of unity.

Let $E$ be an elliptic curve over $k$,
and let $\rho:\Ga_k\to \GL_2(\hat\Z)$ be the Galois
representation in the torsion subgroup of $E$.
The group $\SL_2(\Z/2)$ is isomorphic to the symmetric group $S_3$.
Let us define $\SL_2^+(\hat\Z)$ as the kernel of the composition of
the reduction modulo 2 map $\SL_2(\hat\Z)\to\GL_2(\Z/2)$
with the unique non-trivial homomorphism $\varepsilon:\GL_2(\Z/2)\to\{\pm1\}$.
(The finite group $\GL_2(\Z/2)$ can be identified with the symmetric group
$S_3$, and then $\varepsilon$ is the signature character.)

In this section we prove the following theorem using methods
that go back to Bashmakov \cite[Ch. 5]{Bash}.

\bthe \label{t}
Let $E$ be an elliptic curve over a field $k$ of characteristic zero 
such that $\SL_2^+(\hat\Z)$ is a subgroup of $\rho(\Ga_k)\subset \GL_2(\hat\Z)$.
Let $K$ be a field such that $k\subset K\subset k^{\mathrm{cyc}}$, and let 
$\varphi:E'\to E\times_kK$ 
be an isogeny of elliptic curves over $K$. 
Then for any point $P\in E(K)$ 
that cannot be written as $P=mQ$ with $m>1$ and $Q\in E(K)$, 
the scheme $\varphi^{-1}(P)$ is integral.
\ethe

In the assumption of the theorem, Lemma \ref{p1} (b) below
shows that the isogeny $\varphi:E'\to E$ can be 
identified with the multiplication by $n$ map $E\to E$ 
for some integer $n$. The theorem then follows from Proposition \ref{pp}.

Serre \cite[Prop. 22]{Serre72} proved that for any elliptic curve $E$ 
over $\Q$ the image $\rho(\Ga_\Q)$ is contained in a certain subgroup $H_\Delta$
which only depends on the discriminant $\Delta$ of $E$. The group
$H_\Delta$ has index 2 in $\GL_2(\hat\Z)$ and contains $\SL_2^+(\hat\Z)$.
N. Jones \cite{jones} showed that almost all elliptic curves over $\Q$
are {\em Serre curves} which means that $\rho(\Ga_\Q)=H_\Delta$.
Theorem \ref{t} thus applies to almost all elliptic curves $E$ over $\Q$.

\subsection{Group cohomology}

For a  integer $n>1$ we define $\SL_2^+(\Z/n)$
as $\SL_2(\Z/n)$ if $n$ is odd, and as the kernel of the
following composite map if $n$ is even:
$$\SL_2(\Z/n)\to\SL_2(\Z/2)\to\Z/2,$$
where the first arrow is reduction modulo $2$, and the second arrow
is the signature $\SL_2(\Z/2)\simeq S_3\to\{\pm1\}$.

\bpr \label{p0}
Let $n$ be a positive integer, and let $G\subset\GL_{2}(\Z/n)$
be a subgroup containing $\SL_{2}^+(\Z/n)$.
If $n$ is odd, then $H^1(G, (\Z/n)^2)=0$.
For $n=2^rm$, where $m$ is odd and $r\geq 1$, the abelian group $H^1(G, (\Z/n)^2)$
is annihilated by $2^{r-1}$.
\epr

The proof of Proposition \ref{p0} is based on a few lemmas.

\ble\label{fixp} For any integer $n>1$ we have
$H^0(\SL_{2}^+(\Z/n), (\Z/n)^2)=0$.
\ele
{\em Proof.} 
The group $\SL_{2}^+(\Z/n)$ contains the transformation $(x,y) \mapsto (x+y,-x)$.
$\Box$

\ble \label{odin}
Let  $G_1$ and $G_2$ be finite groups, let $M_1$ be a $G_1$-module
and let $M_2$ be a $G_2$-module such that $(M_1)^{G_1}=(M_2)^{G_2}=0$.
The following natural map is injective:
$$H^1(G_1\times G_2,M_1 \oplus  M_2)\lra
H^1(G_1,M_1)\oplus H^1(G_2,M_2).$$
\ele
{\em Proof.} 
It is enough to prove that $H^1(G_1\times G_2,M_i)$ injects into 
$H^1(G_i,M_i)$ for $i=1,2$. To fix ideas, assume $i=1$. 
We have an inflation-restriction exact sequence
$$0\to H^1(G_2,(M_1)^{G_1}) \to  H^1(G_1\times  G_2, M_1)\to H^1(G_1, M_1),$$
which implies the lemma. $\Box$

\ble \label{dva}
For an odd prime $p$ we have $H^1(\SL_{2}(\Z/p^r), (\Z/p^r)^2)=0$
for any positive integer $r$.
\ele
{\em Proof.} 
Let $C=\{\pm\Id\}\subset\GL_{2}(\Z/p^r)$. 
This is a central subgroup, so we
have an inflation-restriction exact sequence
$$0\to H^1(\SL_{2}(\Z/p^r)/C,((\Z/p^r)^2)^C) \to  
H^1(\SL_{2}(\Z/p^r),(\Z/p^r)^2) \to  H^1(C,(\Z/p^r)^2).$$
The order of $C$ is $2$ but $p$ is odd, so we have $H^1(C,(\Z/p^r)^2)=0$.
We also have $((\Z/p^r)^2)^C=0$, and the lemma follows. $\Box$

\ble \label{tri}
For a positive integer $r$ the group $H^1(\SL_2^+(\Z/2^r), (\Z/2^r)^2)$
is annihilated by $2^{r-1}$.
\ele
{\em Proof.} When $r=1$ the group $\SL_2^+(\Z/2^r)$ has order 3,
so the claim is obvious.

Now suppose $r\geq 2$.
We denote the tautological $\SL_{2}^+(\Z/2^r)$-module 
$(\Z/2^r)^2$ by $M$. Let $\sigma$ be the scalar $2\times 2$-matrix 
$(1+2^{r-1})\Id$, and let $H=\{\Id,\sigma\}$.
It is clear that $H\simeq \Z/2$ is a central subgroup of 
$\SL_{2}^+(\Z/2^r)$ and $M^H=2M$. Let $G=\SL_{2}^+(\Z/2^r)/H$.
There is an inflation-restriction exact sequence 
$$0\to H^1(G,2M) \to H^1(\SL_{2}^+(\Z/2^r),M) \to H^1(H,M)^G.$$
Since $2^{r-1}(2M)=0$ the group
$H^1(G,2M)$ is annihilated by $2^{r-1}$. We have 
$$H^1(H,M)=\Ker[(1+\sigma):M\to M]/(1-\sigma)M.$$
For $r\geq 3$ we have $1+2^{r-2}\in(\Z/2^r)^*$, hence
the kernel of $1+\sigma=2(1+2^{r-2})\Id$ 
is $2^{r-1}M=(1-\sigma)M$. Thus $H^1(H,M)=0$ for $r\geq 3$.
For $r=2$ the map $(1+\sigma):M\to M$ is the multiplication
by $4$ on $M=(\Z/4)^2$, hence in this case $H^1(H,M)=M/2M$.
Since $H$ is central in $\SL_{2}^+(\Z/2^r)$, the action of $G$ on $H$
is trivial, hence $G$ acts on $M/2M$ through its quotient $\SL_{2}^+(\Z/2)$.
The only invariant element under this action is zero.
Thus $H^1(H,M)^G=0$ in all cases, so the lemma is proved. $\Box$

\medskip  

\noindent{\em Proof of Proposition \ref{p0}}. If $R_1$ and $R_2$ 
are commutative rings with $1$, then we have
$\SL_2(R_1\times R_2)\cong \SL_2(R_1)\times \SL_2(R_2)$.
By Lemma \ref{fixp} we have  $H^0(\SL_2^+(\Z/m), (\Z/m)^2)=0$ 
for every positive integer $m$.
Writing $n$ as a product of prime powers,
and applying Lemmas \ref{odin}, \ref{dva} and \ref{tri} we prove
Proposition \ref{p0} in the case $G=\SL_{2}^+(\Z/n)$.
In the general case we note that $\SL_{2}^+(\Z/n)$ is normal in $G$, and
the only $\SL_{2}^+(\Z/n)$-invariant element in $(\Z/n)^2$ is zero.
Hence the restriction from $G$ to $\SL_{2}^+(\Z/n)$ gives an
injective map
$$H^1(G,(\Z/n)^2)\hookrightarrow H^1(\SL_{2}^+(\Z/n),(\Z/n)^2),$$
and the proposition follows.
$\Box$

\subsection{Isogenies of elliptic curves}

The multiplication by $n$ map on an elliptic curve $E$ is
denoted by $[n]:E\to E$. Let $\rho_n:\Ga_k\to \GL_2(\Z/n)$ be the Galois
representation in the $n$-torsion subgroup $E[n]$. 

\ble \label{p1}
Let $E$ be an elliptic curve
over a field $k$ of characteristic zero such that 
$\SL_2^+(\Z/\ell)\subset\rho_\ell(\Ga_k)$ for every prime $\ell$.
Then we have the following statements.

$(i)$ Let $M$ be a  finite $\Ga_k$-submodule of $E(\bar k)$.
Then $M=E[m]$ for an
integer $m\not=0$. In particular $E(k)$ is torsion-free.

$(ii)$ Let $\varphi:C\to E$ be an isogeny of elliptic curves.
Then there is an integer $n>0$ 
and an isomorphism of elliptic curves $\psi:E\tilde\lra C$ such that 
$\varphi\circ\psi=[n]$.
\ele
{\em Proof.} $(i)$ Let $n$ be the smallest positive integer
such that $[n]M=0$. We claim that $M=E[n]$. The $\Ga_k$-module $M$
is the direct sum of its $\ell$-primary  components $M\{\ell\}$.
If $\ell^r$ is the highest power of the prime number $\ell$
that divides $n$, then $\ell^r$ is the smallest positive integer that
annihilates $M\{\ell\}$. The tautological $\SL_2^+(\Z/\ell)$-module 
$(\Z/\ell)^2$ is simple, hence
the $\Ga_k$-module $E[\ell]$ is simple by assumption. But
$\ell^{r-1}M\{\ell\}\not=0$, so we have $\ell^{r-1}M\{\ell\}=E[\ell]$. 
The abelian
group $E[\ell^r]$ is isomorphic to $(\Z/\ell^r\Z)^2$, so any subgroup
that is mapped by multiplication by $\ell^{r-1}$ onto
$(\ell^{r-1}\Z/\ell^r\Z)^2$ is equal to the whole group. Thus $M\{\ell\}=E[\ell^r]$,
and hence $M=E[n]$.

$(ii)$ Passing to the dual isogeny of $\varphi:C\to E$
we see that it is enough to prove that 
every isogeny $\alpha:E\to E'$ is $[m]:E\to E$ for some integer $m$.
This follows from $(i)$ by taking $M=\Ker(\alpha)$. $\Box$

\bpr \label{pp}
Let $E$ be an elliptic curve over a field $k$ of characteristic zero 
such that $\rho(\Ga_k)\subset \GL_2(\hat\Z)$ contains $\SL_2^+(\hat\Z)$.
For any field $K$ such that $k\subset K\subset k^{\mathrm{cyc}}$ 
we have the following statements.

$(i)$ $\SL_2^+(\hat\Z)\subset\rho(\Ga_K)$, hence 
$\SL_2^+(\Z/n)\subset\rho_n(\Ga_K)$ for any positive integer $n$.

$(ii)$ 
For any point $P\in E(K)$ that cannot be written as 
$P=mQ$ for $m>1$ with $Q\in E(K)$, and for any integer $n>0$, 
the scheme $[n]^{-1}(P)$ is integral.
\epr
{\em Proof.} $(i)$ The composition
$\det\rho:\Ga_k\to\GL_{2}(\hat\Z)\to \hat\Z^*$ is the cyclotomic character,
so $\Ga_{k^{\mathrm{cyc}}}$ is the subgroup of $\Ga_k$ given by 
the condition $\det\rho(x)=1$. In view of the natural surjections
$\SL_2(\hat\Z)\to \SL_{2}(\Z/n)$, this proves $(i)$.

$(ii)$ Let $K_n=K(E[n])$. In $(i)$ we proved that
$G=\Ga_K/\Ga_{K_n}$ is a subgroup of $\GL_2(\Z/n)$
containing $\SL_2^+(\Z/n)$. 
The scheme $[n]^{-1}(P)$ is a $K$-torsor of $E[n]$ of the
class $\kappa(P)\in H^1(K,E[n])$, where 
$\kappa:E(K)/n\hookrightarrow  H^1(K,E[n])$ is the Kummer map.
There is an inflation-restriction exact sequence
\begin{equation}
0\to H^1(G,E[n])\to H^1(K,E[n])\to H^1(K_n,E[n])^G.
\label{s1}
\end{equation}
The restriction of $\kappa(P)$ to $H^1(K_n,E[n])$ is 
a $\Ga_K$-equivariant homomorphism $\varphi: \Ga_{K_n}\to E[n]$
(where $\Ga_K$ acts on $\Ga_{K_n}$ by conjugations, and on $E[n]$
in the usual way).
Its image is thus  a $\Ga_K$-submodule of  $E[n]$.
By $(i)$ and by Lemma \ref{p1} $(i)$,
we have $\varphi(\Ga_{K_n})=E[m]$ for some $m|n$.
The set of $\ov K$-points of $[n]^{-1}(P)\times_K K_n$ with
a natural action of $\Ga_{K_n}$ can be identified, by a choice of the
base point, with the set $E[n]$ on which $g\in \Ga_{K_n}$ acts by 
translation by $\varphi(g)$. 

Write $n=2^r s$, where $r\geq 0$ and $s$ is odd. We first deal with
the case $r=0$, and then proceed by induction in $r$. 
Since $P$ is not divisible in $E(K)$ by assumption, 
and $E(K)$ is torsion-free, the order of $\kappa(P) \in H^1(K,E[n])$ is $n$.

If $n$ is odd, by Proposition \ref{p0} the exact sequence (\ref{s1})
gives rise to an embedding of $H^1(K,E[n])$ into
$\Hom(\Ga_{K_n},E[n])$, so that $\varphi$ has order $n$.
It follows that $\varphi(\Ga_{K_n})=E[n]$ 
in this case. This implies that $[n]^{-1}(P)\times_K K_n$ is irreducible,
hence $[n]^{-1}(P)$ is irreducible.

Now suppose that $n=2n'$ and the scheme $[n']^{-1}(P)$ is irreducible.
The multiplication by $2$ map defines a surjective morphism
$[n]^{-1}(P)\to[n']^{-1}(P)$ which is a torsor of $E[2]$.
We know that $\Ga_K$ acts transitively on $[n']^{-1}(P)(\ov K)$
and we want to show that $\Ga_K$ acts transitively on $[n]^{-1}(P)(\ov K)$.
For this we must show that each $\ov K$-fibre of $[n]^{-1}(P)\to[n']^{-1}(P)$
is contained in one $\Ga_K$-orbit. Recall that
the $\Ga_{K_n}$-set $[n]^{-1}(P)(\ov K)$ is identified
with $E[n]$ so that $g\in \Ga_{K_n}$ acts as the translation
by $\varphi(g)$. The $\ov K$-fibres of $[n]^{-1}(P)\to[n']^{-1}(P)$
are the $E[2]$-orbits in $E[n]$. Therefore,
it is enough to show that $\varphi(\Ga_{K_n})$ contains $E[2]$.
As the order of $\kappa(P)$ is $n$, by Proposition \ref{p0} 
the exact sequence (\ref{s1}) shows that the order of $\varphi$
is divisible by $2s$. Thus $\varphi(\Ga_{K_n})$ contains $E[2s]$ 
and hence contains $E[2]$. This finishes the proof. $\Box$

\section{Surfaces} \label{surfaces}

\subsection{An elliptic curve}

For an elliptic curve $E$ over a field $k$ of characteristic zero
we denote by $E^c$ the quadratic twist of $E$ by $c\in k^*$.

\ble \label{le1}
Let $E$ be an elliptic curve over a field $k$ of characteristic zero. 
For a quadratic extension $K=k(\sqrt{d})$ we have an exact sequence
$$0\to E^d(k)\to E(K)\to E(k).$$
\ele
{\em Proof.} Let $\sigma$ be the non-zero element of $\Gal(K/k)$.
We have $E(k)\simeq E(K)^\sigma$.
The choice of a square root of $d$ in $K$ defines 
an isomorphism $E^d\times_k K\simeq E\times_k K$. This gives an
identification $E^d(k)\simeq\{x\in E(K)|\sigma(x)=-x\}$.
Sending $x\in E(K)$ to $x+\sigma(x)$ defines a homomorphism
$E(K)\to E(k)$ with kernel $E^d(k)$. $\Box$

\bpr \label{p2}
There exist the following data:

a real quadratic field $k=\Q(\sqrt{c})$, where $c$ 
is a square-free positive integer not congruent to $1$ modulo $8$; 

a totally real biquadratic field $K=\Q(\sqrt{c},\sqrt{d})$, 
where $d$ is a square-free positive integer;

an elliptic curve $E$ over $\Q$ of discriminant
$\Delta<0$, such that $\SL_2^+(\hat\Z)\subset \rho(\Ga_k)$,
$E(k)=\{0\}$, and  $E(K)$ is torsion-free of positive rank.
\epr
{\em Proof.} Let $E$ be the curve $y^2+y=x^3+x^2-12x-21$ of conductor $67$
and discriminant $\Delta=-67$, and take
$c=10$, $d=2$. Using {\tt sage} we check that 
$E(\Q)=E^{10}(\Q)=0$. By Lemma \ref{le1} we have $E(k)=0$.
Using {\tt sage} we check that $E^{2}(\Q)\simeq E^{5}(\Q)\simeq\Z$, so
by Lemma \ref{le1} we conclude that $E(K)$ is torsion-free of 
positive rank. We claim that $E$ is a Serre curve, which means that
$\rho(\Ga_\Q)=H_\Delta$, where $H_\Delta$ is a subgroup of $\GL_2(\hat\Z)$
of index 2 containing $\SL_2^+(\hat\Z)$. By a result of N. Jones 
\cite[Lemma 5]{jones}, for this it is enough to show that 
$E$ satisfies the following conditions: 

$\rho_\ell(\Ga_\Q)=\GL_2(\Z/\ell)$ for all primes $\ell$
(this is checked using {\tt sage}), 

$\rho_8(\Ga_\Q)=\GL_2(\Z/8)$ (this is checked using \cite{Dok}), and 

$\rho_9(\Ga_\Q)=\GL_2(\Z/9)$ (this is checked using \cite{Elk}). 

\noindent 
By Proposition \ref{pp} $(i)$, $\SL_2^+(\hat\Z) \subset \rho(\Ga_\Q)$
implies $\SL_2^+(\hat\Z) \subset \rho(\Ga_k)$.
This finishes the proof. $\Box$

\begin{rem} \label{r1}
For $c$ and  $k=\Q(\sqrt{c})$ as above,
the conic $x^2+y^2+z^2=0$ has points in all non-Archimedean
completions of $k$, but no points in the two real completions of $k$. 
Indeed, this conic has $\Q_p$-points
for all odd primes $p$. Since 2 is ramified or inert in $k$,
this conic also has points in the completion of $k$ at the unique
prime over $2$.
\end{rem}

\subsection{A conic bundle over the elliptic curve} \label{dim2}
 
Let $k$, $K$ and $E$ be as in Proposition \ref{p2}.
The elliptic curve $E$ can be given by its short Weierstra{\ss} equation
$$y^2=r(x),$$
where $r(x)=x^3+px+q$, for $p,q \in \Q$. 
Since $\Delta=-4p^3-27q^2<0$ the topological space $E(\R)$ is connected.
The neutral element of $E(k)$ is the point at infinity.
We denote by $\pi:E\to\P^1_k$ the projection sending $(x,y)$ to $x$.

Let $P\in E(K)$ be a point not divisible in $E(K)$.
Let $\sigma \in \Gal(K/k)$ be the generator.
Since $P+\sigma(P)\in E(k)=0$ we obtain that
$\pi(P)$ is a point of $\A^1_k(k)=k$,
say $\pi(P)=a \in k$. The $K$-point $P$ gives rise to 
a solution of $y^2=r(a)$ in $K$, hence
$r(a) \in k$ is totally positive. We have $r(a)\not=0$ since 
$E(K)$ is torsion-free.

Let $b \in k$, $b \neq a$.   
We define a central simple algebra $A$ over $k(\P_k^1)=k(x)$ 
as a tensor product of quaternion algebras
\begin{equation}
A= ((x-a)/(x-b),r(b))\otimes(-1,-1).\label{A}
\end{equation}
The algebra $A$ is unramified outside of the points $x=a$ and $x=b$.
At each of these points the residue of $A$ is given by the class of 
$r(b)$ in $k^*/k^{*2}$. 

\bpr[Albert]  \label{Albert} 
Let $F$ be a field, $char(F)\neq 2$, and let 
$\alpha,\beta,\gamma,\delta  \in F^*$. The tensor product of quaternion algebras
$(\alpha,\beta) \otimes (\gamma,\delta )$ is a division algebra if and 
only if the diagonal quadratic form 
$\langle\alpha,\beta,-\alpha\beta,-\gamma,-\delta ,\gamma\delta \rangle$ 
is anisotropic. 
If it is isotropic, then $(\alpha,\beta) \otimes (\gamma,\delta )$ 
is similar to a quaternion algebra over $F$.
\epr

For the proof see \cite[\S 16.A]{KMRT}.
The quadratic form 
$\langle \alpha,\beta,-\alpha\beta,-\gamma,-\delta ,\gamma\delta \rangle$
is called an {\em Albert form} associated to 
$(\alpha,\beta) \otimes (\gamma,\delta )$.

\ble \label{le2}
Let $k$, $K$, $a$, $r(t)$  be as above.
If $b \in k$ is such that $r(b)$ is totally negative, then
the algebra  $A$ over the field 
$k(\P_k^1)=k(x)$ is similar to a quaternion algebra.
\ele
{\em Proof.} The associated Albert form contains the subform
$\Phi=\langle r(b),1,1,1\rangle$. By Remark \ref{r1}
the form $\langle 1,1,1\rangle$ 
is isotropic over all finite completions of $k$.
Since $r(b)$ is totally negative, $\Phi$
is isotropic over both real completions of $k$.
By the Hasse--Minkowski theorem the quadratic form 
$\Phi$ is isotropic over $k$. 
An application of Proposition \ref{Albert} concludes the proof. $\Box$

\medskip

We are now ready to state one of the main results of this paper.
(For an explanation why we do not consider the case $k=\Q$
see Proposition \ref{last} below.)

\bthe \label{t1}
There exist a real quadratic field $k$, an elliptic curve $E$ 
and a smooth, projective and geometrically integral surface $X$ over $k$
with a surjective morphism $f :X\to E$ satisfying the following properties:

$(i)$ the fibres of $f :X\to E$ are conics;

$(ii)$ there exists a closed point $P\in E$ such that the field  $k(P)$
is a totally real biquadratic extension of $\Q$
and the restriction $X\setminus f^{-1}(P)\to E\setminus P$ is a smooth morphism;

$(iii)$ $X(\AA_k)^{\et,\Br}\not=\emptyset$ and $X(k)=\emptyset$.
\ethe
{\em Proof.} 
We keep the above notation: $k=\Q(\sqrt{10})$,
$E$ is the curve of Proposition \ref{p2} given by its 
short Weierstra{\ss} equation
$y^2=r(x)$, viewed as a curve over $k$,
and $K=k(P)=\Q(\sqrt{10},\sqrt{2})$.
We fix $b \in k$ such that $r(b)$ is totally negative, and
define a central simple algebra $A$ by (\ref{A}).
By Lemma \ref{le2}, $A$ is similar to a quaternion algebra $B$
over $k(\P_k^1)$. Let $f :X \to E$ be a 
relatively minimal conic bundle 
such that the generic fibre $X_{k(E)}$ is the conic 
over $k(E)$ defined by the quaternion algebra $B\otimes_{k(\P^1)}k(E)$.
(We call a conic bundle $X\to E$ relatively minimal if 
for every conic bundle $Y\to E$ any birational morphism $X\to Y$
compatible with the projections to $E$, is an isomorphism.
See \cite[Thm. 1.6]{Manin66} for the well known description
of the fibres of $X\to E$.)
The closed point $P=\pi^{-1}(a)\simeq\Spec\,K$ is
a solution to $y^2=r(a)$. The residue of $B\otimes_{k(\P^1)}k(E)$ at $P$
is the class of $r(b)$ in $K^{\times}/K^{\times 2}$, 
which is non-trivial since $r(b)$ is totally negative. We have 
$\pi^{-1}(b)=\Spec\,k(\sqrt{r(b)})$, and the residue at this 
closed point
is the class of $r(b)$ in $k(\sqrt{r(b)})^*/k(\sqrt{r(b)})^{*2}$ which 
is trivial. Thus $f :X \to E$ is smooth away from $P$.
This gives $(i)$ and $(ii)$.

If we go over to one  of the  two real completions $k_v$ of $k$, 
the point $P$
breaks up into two real points $P_{1}$ and $P_{2}$. The residue
at each of these points is an element of $\R^*/\R^{*2}$ 
given by the image of $r(b)$, and that element is totally negative.
Thus the algebra $A$ is ramified at both $P_1$ and $P_2$.
The fibre of $X \to E$ above each point $P_{i} \in E(k_{v})$
is thus a singular conic, hence has a $k_{v}$-point. In particular,
$X(k_v) \neq \emptyset$ for each of the real completions of $k$. 
The condition $\Delta<0$ implies that
$E(k_v)$ is connected, thus any point in the
image of $X(k_v) \to E(k_v)$ is path connected to
the point $0$, which is the unique point of $E(k)$. (In fact,
the image of $X(\R)$ in $E(\R)$ is the closed interval between $P_{1}$
and $P_{2} $ which does not contain $0$.)
Let $M_v$ be any point of $X(k_v)$, and let
$I_v\subset E(k_v)\simeq S^1$ be a real interval linking
$f(M_v)$ and $0$.

The value of $A$ at $\infty\in\P^1_{k} $, hence also at $0 \in E(k)$,
is $(-1,-1)$, hence the fibre
$X_0=f^{-1}(0)$ is the conic $x^2+y^2+z^2=0$. By Remark \ref{r1}, $X_0$ has
points in all completions of $k$ except the  two real
completions, hence $X_0(k)=\emptyset$. Since $E(k)=\{0\}$
it follows that $X(\AA_k)\not=\emptyset$, but $X(k)=\emptyset$.
For each finite place $v$ of $k$ choose $M_v\in X_0(k_v)$.

We now prove that $(M_v)\in X(\AA_k)^{\et,\Br}$.
By Proposition \ref{et-cov} every torsor $X'\to X$ of a finite  
$k$-group scheme $G$ is the pullback 
$X'=X\times_E E\to X$ of a torsor $E'\to E$ of $G$.
By twisting $E'$ and $X'$ with a $k$-torsor of $G$ (and replacing $G$
with the corresponding inner form) we can assume
that $E'$ has a $k$-point $0'$ over $0\in E(k)$. The connected
component $C$ of $E'$ containing this $k$-point is a smooth,
projective, geometrically integral curve. The $k$-morphism 
$\varphi:C\to E$ is finite and \'etale, hence $C$ has genus~1.
Choosing $0'$ for the origin of the group law on $C$ 
we make $\varphi:C\to E$ into an isogeny of elliptic curves and
write $0=0'$.

Let $Y=X\times_E C$. Then $Y$ is a smooth, projective, geometrically
integral surface over $k$ which is an irreducible component of $X'$.
The morphism $g :Y\to C$ is a conic bundle.
By Proposition \ref{p2} (3) we have 
$\SL_2^+(\hat\Z)\subset \rho(\Ga_k)$.
Since the point $P \in E(K)$ is not divisible, we can apply Theorem \ref{t}
to the elliptic curve $E$ and the isogeny $\varphi:C\to E $.
It follows that the inverse image $Q=\varphi^{-1}(P)$ 
of the closed point $P$ of the $k$-curve $E$  is a closed 
point of the $k$-curve $C$. 
We see that $g:Y\to C$ is a conic bundle
such that $Y\setminus Y_{Q}\to C\setminus Q$ is a smooth morphism. 
By Proposition \ref{Br-const} $(i)$, the induced map
$g^*:\Br(C)\to\Br(Y)$ is surjective.

For each finite place $v$ of $k$ let $M'_v$ be the $k_v$-point
in the fibre $Y_0$ over $0\in C(k)$ that projects to $M_v\in X_0(k_v)$.
Now let $v$ be a real place of $k$.
By Lemma \ref{p1} $(ii)$ 
the isogeny $\varphi:C\to E$ is identified
with $[n]:E\to E$ for some $n$.
Since $\Delta<0$, the induced map $C(k_{v})\simeq E(k_v)\to E(k_{v})$ 
is a surjective \'etale map $S^1 \to S^1$. Since $I_v$
is contractible, we see that $\varphi^{-1}(I_v)$ is a disjoint union
of copies of $I_v$ exactly one of which contains $0$.
Let us call this interval $I'_v$. One of its ends is $0$ and 
the other end is a point $R\in C(k_v)$ such that $\varphi(R)=f(M_v)$.
Hence the real fibre $Y_R$ is naturally isomorphic to
the fibre of $f :X\to E$ that contains $M_v$. Let $M'_v$
be the point in $Y_R(\R)$ that projects to $M_v$.
Since $\Br(Y)=g^*(\Br(C))$, we see that
$(M'_v)\in Y(\AA_k)^\Br$. But $Y$ is an irreducible component of $X'$
hence $(M'_v)\in X'(\AA_k)^\Br$. This is a lifting of the adelic point $(M_v)$,
and therefore $(M_v)\in X(\AA_k)^{\et,\Br}$. $\Box$

\begin{rem} \label{r3}
In the proof above the ad\`ele $(f(M_v))\in E(\AA_k)$ and the
$k$-point $0 \in E(k)$ have the following property:
for each place $v$ the class of the $0$-cycle 
$f(M_v)-0$ is infinitely divisible in $\Pic(E\times_k k_v)$, in fact
it is even zero if $v$ is non-archimedean.
Methods initiated by one of the authors and developed by 
Frossard \cite[Thm. 0.3]{F}, by van Hamel, and by Wittenberg \cite[Thm. 1.3]{W},
then show that there exists a $0$-cycle of degree $1$ on $X$.
In these various papers,  
finiteness of the Shafarevich--Tate group of 
the Jacobian of the base curve is assumed in order
to appeal to the Cassels--Tate dual exact sequence, which guarantees the
existence of a $0$-cycle of degree $1$ on the curve satisfying 
a divisibility property analogous to the one above.
In our case we can take this $0$-cycle of degree $1$ to be 
the $k$-point $0$, so there is no need to assume the finiteness of the
Shafarevich--Tate group of $E$.
\end{rem}

\section{Remarks on the Brauer--Manin set} \label{varia}

\subsection{Birational invariance}

Recall that $\Br_{0}(X)$ denotes the image of the map $\Br(k) \to \Br(X)$.
 
\bpr \label{5.1}
Let $k$ be a number field, and let $X$ and $Y$ be smooth, projective,
geometrically integral varieties over $k$ that are birationally equivalent.

$(i)$ If $X(\AA_k)^\Br\not=\emptyset$, then $Y(\AA_{k})^\Br\not=\emptyset$.

$(ii)$ If $X(\AA_k)^{\et,\Br}\not=\emptyset$, then $Y(\AA_{k})^{\et,\Br}\not=\emptyset$.

$(iii)$ Assume, in addition, that $\Br(X)/\Br_0(X)$ is finite. Then the
density of $X(k)$ in $X(\AA_k)^\Br$ implies the density of $Y(k)$ in $Y(\AA_k)^\Br$.
\epr
{\em Proof.} $(i)$ By Hironaka's theorem, there exists a smooth, 
projective variety $Z$ over $k$ and birational morphisms 
$f : Z \to X$ and $g: Z \to Y$. Let $v$ be a place of $k$.
Since $Z$ is proper, $f(Z(k_{v}))$ is closed in $X(k_{v})$. 
Let us show that $X(k_{v})=f(Z(k_{v}))$, for which it is enough 
to show that $f(Z(k_{v}))$ is dense in $X(k_{v})$.
There exists a non-empty Zariski open set $U\subset X$ 
such that $f$ induces an isomorphism $f^{-1}(U)\tilde\lra U$.
Then $U(k_{v})\subset f(Z(k_{v}))$, but $U(k_{v})$ is dense in $X(k_{v})$
by the implicit function theorem. 

Thus for any $(M_{v}) \in X(\AA_k)$ there exists 
$(N_{v}) \in Z(\AA_k)$ such that $(f(M_v))=(N_v)$.
The birational morphism $f$ induces an isomorphism
$f^* : \Br(X) \tilde\lra \Br(Z)$. From the projection formula
we conclude that if $(M_{v})\in X(\AA_{k})^{\Br}$, then 
$(N_{v})\in Z(\AA_{k})^{\Br}$. By the covariant functoriality of the
Brauer--Manin set we have $g\left(Z(\AA_{k})^{\Br}\right)
\subset Y(\AA_{k})^{\Br}$.

$(ii)$ Let $G$ be a finite $k$-group scheme.
By the birational equivalence of the fundamental group
\cite[X, Cor. 3.4]{SGA1} there is a natural bijection
between $X$-torsors and $Y$-torsors of $G$ in which a
torsor $X'/X$ corresponds to $Y'/Y$ if $X'\times_X Z=Y'\times_Y Z$.
(This bijection respects the twisting by a $k$-torsor of $G$.)
Let us denote this $Z$-torsor of $G$ by $Z'$. 
The natural morphism $Z'\to X'$ is
a componentwise birational morphism of
smooth and projective varieties,
so it induces an isomorphism $\Br(X') \tilde\lra \Br(Z')$.

Consider any $(M_v) \in X(\AA_k)^{\et,\Br}$. 
Let $(N_v)\in Z(\AA_k)$ be a lifting of $(M_v)$ as in part $(i)$.
By the definition of $X(\AA_k)^{\et,\Br}$,
for any torsor $X'/X$ the ad\`ele $(M_v)$ is the image of some 
$(M'_v) \in X'(\AA_k)^\Br$ (after twisting $X'/X$ by
a $k$-torsor of $G$ and replacing $G$ by the corresponding inner form). 
If $N'_v\in Z'(k_v)$ 
is the point $M'_v\times_{M_v}N_v$ for each place $v$, then
$(N'_v)\in Z'(\AA_k)^\Br$ by the projection 
formula. This implies that $(N_v)\in Z(\AA_k)^{\et,\Br}$.
By the covariant functoriality of the \'etale Brauer--Manin set
we have $g(Z(\AA_k)^{\et,\Br})\subset Y(\AA_k)^{\et,\Br}$.

$(iii)$ Let $U \subset X$ and $V \subset Y$ be non-empty Zariski open sets
such that there is an isomorphism $h : V \tilde\lra U$.
Under the assumption of $(ii)$, the group $\Br(Y)/\Br_0(Y) \simeq \Br(X)/\Br_0(X)$
is finite. Let $\beta_{1}, \dots, \beta_{n} \in \Br(Y)$ be 
coset representatives. There exists a finite set $\Sigma$ of places of 
$k$ such that each $\beta_{i}$ vanishes on $Y(k_{v})$ for each 
$v\notin \Sigma$. Thus every ad\`ele $(M_v)\in Y(\AA_k)^\Br$
can be approximated by an ad\`ele $(N_v)\in Y(\AA_k)^\Br$ 
such that $N_v\in V(k_v)$ for all completions $k_v$ of $k$.
Let $M_{v}=h(N_{v})$. Since $h$ induces an isomorphism $\Br(Y)\simeq\Br(X)$,
the ad\`ele $(M_{v})$ is in $X(\AA_k)^\Br$. If $X(k)$ is dense 
in $X(\AA_k)^\Br$, we can approximate $(M_{v})$
by a $k$-point $M$ in $U$. Then $N=h^{-1}(M) \in V(k)$
is close to $(N_v)\in Y(\AA_k)^\Br$.  $\Box$

\begin{rem} (1) We do not know if the analogue of 
Proposition \ref{5.1} $(iii)$ holds for the \'etale Brauer--Manin set.

(2) We do not know if Proposition \ref{5.1} $(iii)$ still holds when
$\Br(X)/\Br_0(X)$ is infinite,
but we can make the following observation.
Recall that $X(\AA_{k})_{\bullet}$ denotes the quotient of $X(\AA_{k})$
by the relation which identifies two points in the same connected component.
{\it For smooth, projective, geometrically integral varieties $X$ and $Y$
that are birationally equivalent,
if one does not assume the finiteness of $\Br(X)/\Br_0(X)$,
then $Y(k)$ can be dense in $Y(\AA_k)^\Br_\bullet$ without
$X(k)$ being dense in $X(\AA_k)^\Br_\bullet$.}
Indeed, let $Y=A$ be an abelian surface over $\Q$ such that $A(\Q)=
\{0\}$ and the Shafarevich--Tate group of $A$ is finite.
Let $Z$ be the blowing-up of $A$ at the $\Q$-point $0$,
and let $X$ be the blowing-up of $Z$ at some $\Q$-point
(all $\Q$-points of $Z$ are contained in the exceptional divisor).
The surface $X$ contains two copies of $\P^1_{\Q}$ meeting at a
$\Q$-point $P$; let us call them $E$ and $F$.
The finiteness of the Shafarevich--Tate group of $A$ implies
\cite[Prop. 6.2.4]{Sk} that $A(\AA_\Q)^\Br$ is the 
connected component of $0$ in $A(\AA_\Q)$, which is isomorphic
to the real connected component of $0$ in $A(\R)$.
Choose a $\Q$-point $M\not=P$ in $E$, and a $\Q$-point 
$N\not=P$ in $F$. Let $q$ be a prime.
Consider the adelic point $(M_p)$ of $X$, where
$M_p=M$ for all $p\not=q$ (including $p=\infty$) and $M_q=N$.
Since the morphism $f:X\to A$ induces an isomorphism
$f^*:\Br(A)\tilde\lra\Br(X)$ and $f(M_p)=0$ for all $p$,
we see that $(M_p) \in X(\AA_\Q)^\Br$.
However, the connected component of $(M_p)$ 
in $X(\AA_\Q)$ contains no $\Q$-points. Indeed,
a $\Q$-point of $X$ is either in $E$ or in $F$.
In the first case it cannot approximate $ M_q$
in the $q$-adic topology, and in the second case it
cannot approximate $ M_p$ in the $p$-adic
topology where $p$ is any prime different from $q$.
It is easy to see that $(M_p) \in X(\AA_\Q)^{\et,\Br}$,
so previous discussion applies to the set of connected components
of the \'etale Brauer--Manin set as well.
\end{rem}

\subsection{Cases where the Brauer--Manin obstruction suffices}

The following result shows that for the conclusion of Theorem \ref{t1}
to hold, the conic bundle $f : X \to E$ must contain degenerate fibres.
 
\bpr Let $E$ be an elliptic curve over a number field $k$
with a finite Shafarevich--Tate group.
Let $f : X \to E$ be a Severi--Brauer scheme over $E$. 
Then $X(\AA_{k})^{\Br} \neq \emptyset$ implies $X(k) \neq \emptyset$.
Moreover, $X(k)$ is dense in $X(\AA_{k})_{\bullet}^{\Br}$.
\epr
{\em Proof.} Since $f : X \to E$ is a projective morphism with 
smooth geometrically integral fibres, there exists a finite
set of places $\Sigma$ such that $E(k_v)=f(X(k_v))$ for $v\notin\Sigma$.
We may assume that $\Sigma$ contains the archimedean places of $k$.
At an arbitrary place $v$ the set $f(X(k_{v}))$ is open and closed in
$E(k_{v})$. Let $(M_{v}) \in X(\AA_{k})^{\Br}$. By functoriality
we then have $(f(M_{v}))\in E(\AA_{k})^{\Br}$.  
The finiteness of the Shafarevich--Tate group of $E$
implies  \cite[Prop. 6.2.4]{Sk} the
exactness of the Cassels--Tate dual sequence
\begin{equation}
0 \to E(k)\otimes \hat{\Z} \to \prod E(k_{v})_{\bullet} \to 
\Hom(\Br (E), \Q/\Z), \label{CT}
\end{equation}
where $E(k_{v})_{\bullet}=E(k_{v})$ if $v$ is a finite place of $k$,
and $E(k_{v})_{\bullet}=\pi_0(E(k_v))$ if $v$ is an archimedean
place. By a theorem of Serre, the image of $E(k)\otimes\hat\Z$ 
in that product coincides with the topological closure of $E(k)$,
see \cite{Serre}, \cite{Wang}.
Approximating at the places of $\Sigma$, we find
a $k$-point $M \in E(k)$ such that the fibre $X_{M}=f^{-1}(M)$
is a Severi--Brauer variety with points in all $k_{v}$ 
for $v \in \Sigma$, hence also for all places $v$. Since $X_{M}$ is 
a Severi--Brauer variety over $k$, it contains a $k$-point, 
hence $X(k)\neq \emptyset$. For the last statement of the theorem
we include into $\Sigma$ the places where we want to approximate.
If $k_v\simeq\R$, each connected component $X(k_v)$ maps
surjectively onto a connected component of $E(k_v)$.
The Severi--Brauer varieties satisfy 
the Hasse principle and weak approximation, so an application
of the implicit function theorem finishes the proof. $\Box$
  
\begin{rem}
The same argument works more generally
for a projective morphism $f: X \to E$
such that each fibre contains a geometrically integral component of 
multiplicity $1$, provided that the smooth $k$-fibres satisfy 
the Hasse principle. For the last statement to hold, 
the smooth $k$-fibres also need to satisfy weak approximation.
\end{rem}

The following proposition is a complement to Theorem \ref{t1} 
which explains why a similar counterexample cannot be constructed 
over $\Q$. 

\bpr \label{last}
Let $E$ be an elliptic curve over a number field $k$
such that both $E(k)$ and the Shafarevich--Tate group of $E$
are finite. Let $f : X \to E$ be a conic bundle.
Suppose that there exists a real place $v_0$ of $k$ such that
for each real place $v\not=v_0$ no
singular fibre of $f : X \to E$ is over a $k_v$-point of $E$.
Then $X(\AA_{k})^{\Br}\neq\emptyset$ implies $X(k)\neq\emptyset$.
\epr
{\em Proof.}  If a $k$-fibre of $f$ is not smooth,
then this fibre contains a $k$-point. We may thus assume
that the fibres above $E(k)$ are smooth.
Let $(M_{v})\in X(\AA_{k})^{\Br}$. Then $(f(M_{v}))\in E(\AA_{k})^{\Br}$. 
Set $N_{v}=f(M_{v})$ for each place $v$. The
finiteness of the Shafarevich--Tate group of $E$
implies the exactness of the Cassels--Tate dual sequence (\ref{CT}).
Hence there exists $N \in E(k)$ such that $N=N_{v} $ for each finite 
place $v$ and such that $N$ lies in the same connected
component as $N_{v}$ for $v$ archimedean. The fibre $X_{N}$ is a smooth 
conic with points in all finite completions of $k$. 
For an archimedean place $v\not=v_0$, the map $X(k_{v} ) \to E(k_{v})$ 
sends each connected component of $X(k_{v})$ onto a connected component 
of $E(k_{v})$. Since $N$ and $N_{v}$ are in the same connected 
component of $E(k_{v})$, this implies $X_{N}(k_{v}) \neq \emptyset$.
Thus the conic $X_{N}$ has points in all completions of $k$ 
except possibly $k_{v_0}$. By the reciprocity law
it has points in all completions of $k$ and hence in $k$.  $\Box$

\begin{rem}
If $E(k)$ is finite
we cannot in general expect $X(k)$ to be dense in $X(\AA_k)_\bullet^\Br$
or even in $X(\AA_k)_\bullet^{\et,\Br}$. 
Indeed, if the fibre $X_N$ over some $k$-point $N$ of $E$ is
an {\em irreducible singular} conic, then the singular point $P$ of $X_N$
is the unique $k$-point of $X_N$. Let $(M_v)$ be any ad\`ele in $X_N(\AA_k)$.
Note that if $v$ splits in the quadratic extension of $k$
over which the components of $X_N$ are defined, then
$X_N\times_k k_v$ is a union of two projective lines meeting at $P$. 
Using the fact that $\Br(\P^1_{k_v})=\Br(k_v)$ we see that $(M_v)\in X(\AA_k)^\Br$.
Furthermore, we have $(M_v)\in X(\AA_k)^{\et,\Br}$, cf.
\cite[Remark 2.4]{HaSk}.
On the other hand, $(M_v)$ is not in the closure of
$X(k)$ in $X(\AA_k)_\bullet$ provided $M_v\not=P$ for at least one
finite place $v$ of $k$.
\end{rem}

\end{document}